\documentclass[10pt,amscd,amssymb,epsf]{amsart}

\usepackage{graphicx}              

\begin{document}
\title{An integer valued SU(3) Casson invariant}
\author{Hans U. Boden}
\address{Department of Mathematics \& Statistics, McMaster University, 
Hamilton, ON, L8S 4K1 Canada}
\email{boden@icarus.mcmaster.ca}

\author{ Christopher M. Herald}
\address{Department of Mathematics, University of Nevada, Reno NV 89557}
\email{herald@unr.edu}

\author{Paul A. Kirk}
\address{Department of Mathematics, Indiana University, Bloomington, IN 47405}
\email{pkirk@indiana.edu}

\date{August 11, 2000}
\begin{abstract}
{We define an integer valued invariant of homology spheres using
the methods of $SU(3)$ gauge theory and study  its behavior under 
orientation reversal and connected sum.}
\end{abstract}

\maketitle

\newcommand{\const}{\mbox{const}}
\newcommand{\lto}{\longrightarrow}
\newcommand{\al}{\alpha}
\newcommand{\be}{\beta}
\newcommand{\ga}{\gamma}
\newcommand{\de}{\delta}
\newcommand{\ep}{\epsilon}
 \renewcommand{\th}{\theta}
\newcommand{\la}{\lambda}
\newcommand{\om}{\omega}
\newcommand{\si}{\sigma}
\newcommand{\Ga}{\Gamma}
\newcommand{\Om}{\Omega}
\newcommand{\Si}{\Sigma}

\newcommand{\tlam}{\widetilde \lambda}

\newcommand{\ZZ}{{\mathbb Z}}
\newcommand{\RR}{{\mathbb R}}
\newcommand{\CC}{{\mathbb C}}
\newcommand{\QQ}{{\mathbb Q}}
\newcommand{\HH}{{\mathbb H}}

\newcommand{\cA}{{\mathcal A}}
\newcommand{\cB}{{\mathcal B}}
\newcommand{\cF}{{\mathcal F}}
\newcommand{\cG}{{\mathcal G}}
\newcommand{\cM}{{\mathcal M}}
\newcommand{\hA}{\widehat{A}}
\newcommand{\hB}{\widehat{B}}

\newcommand{\tB}{{\widetilde{\mathcal B}}}
\newcommand{\tM}{{\widetilde{\mathcal M}}}
\newcommand{\tC}{{\widetilde{C}}}

\newcommand{\hh}{{\mathfrak h}}
\newcommand{\hhp}{{{\mathfrak h}^\perp}}
\newcommand{\im}{\operatorname{im}}
\newcommand{\hol}{\operatorname{{\it hol}}}
\newcommand{\tr}{\operatorname{\it tr}}
\newcommand{\Spec}{\operatorname{Spec}}
\newcommand{\crit}{\operatorname{Crit}}
\newcommand{\hess}{\operatorname{Hess}}
\newcommand{\ind}{\operatorname{ind}}

 \newcommand{\opX}{{\bar{X}}}

\newtheorem{defn}{Definition}
\newtheorem{lemma}[defn]{Lemma}
\newtheorem{thm}[defn]{Theorem}
\newtheorem{prop}[defn]{Proposition}
 \newtheorem{cor}[defn]{Corollary}

\section{Introduction}

 Any $SU(n)$
generalization of the Casson invariant for homology 3-spheres $X$
ought to be defined 
as a signed count of conjugacy classes  of irreducible $SU(n)$
representations of $\pi_1 X $. The difficulty is that, just as in the $SU(2)$ case, one must
perturb the space of representations to
make it generic and hence finite, but for $n>2$ the signed count depends on
the perturbation. 

In \cite{taubes}, Taubes introduced a technique for perturbing the flatness
equations and gave a gauge-theoretic interpretation of
  Casson's $SU(2)$ invariant
as a signed count of  gauge orbits of
perturbed flat $SU(2)$ connections on $X$.
In \cite{bh1}, an $SU(3)$ Casson invariant $\la_{SU(3)}$ for homology 3-spheres
$X$ is defined using the perturbation  approach of Taubes. 
For each generic
perturbation $h$, an integer $\la'_{SU(3)}(X,h)$  is defined as a signed
count of gauge orbits of  irreducible,
$h$-perturbed flat SU(3) connections,
\begin{equation} \label{ladefn}
\la'_{SU(3)}(X,h) =
\sum_{[A] \in \cM^*_{h}} (-1)^{SF(\th,A)},
\end{equation}
where $\cM^*_{h}$ is the moduli space of {\it irreducible}, $h$-perturbed flat
$SU(3)$ connections, $A$ is a representative for the gauge orbit $[A]$,
and $SF$ refers to the spectral flow of the odd signature operator $K_A$
acting on $su(3)$-valued forms.

The resulting integer $\la'_{SU(3)}(X,h)$   is not independent
of the choice of the perturbation $h$. To extract a topological invariant,
one needs to define
a correction term. An analysis of the
parametrized moduli space corresponding to a path joining two generic
perturbations suggests that  the correction term
should be a signed sum of the form
\begin{equation}\label{lousy}
\tfrac{1}{2} \sum_{[A] \in \cM^r_{h}}
(-1)^{SF(\theta ,A)}SF_\hhp(A_0, A),\end{equation}
where $\cM^r_{h}$ is the moduli space of {\it reducible}, $h$-perturbed flat
connections  (i.e., with holonomy reducing to
$S(U(2)\times U(1))$)
and $A_0$ is some fixed reducible connection.
Here, the subscript on $SF_\hhp$ indicates that  the operator $K_A$ is 
acting on forms  
with coefficients in  $\hhp,$
the orthogonal complement in $su(3)$ of
the Lie subalgebra  $\hh= s(u(2)\times u(1))$.

In \cite{bh1}, $A_0$ was taken to be $\th$, the trivial
connection.  Unfortunately,
 the quantity $SF_{\hhp}(\theta, A)$  is not gauge invariant; it depends
on the choice of representative $A$ for the gauge equivalence class $[A]$.
In \cite{bh1}, it was shown how to restore gauge invariance
by restricting to small perturbations and adding
the Chern-Simons invariant of a flat connection $\widehat A$ near the 
representative $A$
to obtain the correction term:
$$\la''_{SU(3)}(X,h)=\tfrac{1}{2}\sum_{[A] \in \cM^r_{h}}
(-1)^{SF(\th,A)}(SF_\hhp(\th,A)-4 cs (\hA) +2 ).$$
Then  $\la_{SU(3)}(X)=\la' _{SU(3)}(X,h) +\la'' _{SU(3)}(X,h)$ 
is  independent
of the choice of small perturbation $h$.
In the recent
preprint \cite{CLM}, Cappell, Lee and Miller develop a different technique
for correcting the gauge ambiguity in (\ref{lousy}).

In this paper, we construct a correction term $\tau''(X,h)$ like (\ref{lousy})
but without  any gauge ambiguity.
Setting    $$\tau(X) =  \la'_{SU(3)}(X,h) + \tau''(X,h)$$
gives an integer valued $SU(3)$ Casson invariant
 of homology 3-spheres.
Like $\la_{SU(3)}$, the invariant $\tau$ enjoys properties (i) -- (iii) of
the following
theorem, which is our main result.
\vskip.2in
  \noindent{\bf Theorem \ref{theorem4}.} {\it
The quantity $\tau(X)$  is
an integer valued  invariant
of homology 3-spheres. Furthermore,
\begin{enumerate}
\item[(i)] If $\tau(X) \neq 0, $ then there exists an irreducible
representation $\rho: \pi_1 X \to SU(3).$
\item[(ii)] If $\bar{X}$ equals $X$ with the orientation reversed, 
then $\tau(X) = \tau(\bar{X})$.
\item[(iii)]  If $X_1$ and $X_2$ are homology 3-spheres, then \\
$\tau (X_1 \# X_2) = \tau(X_1) + \tau (X_2) + 4\la_{SU(2)}(X_1)
\la_{SU(2)} (X_2).$
\end{enumerate}}
  \vskip.2in
The invariant $\tau$ has numerous advantages over $\la_{SU(3)}$, and we mention three.
First, $\tau(X) \in \ZZ$.
(A priori  $\la_{SU(3)}(X)$ takes values in $\RR$, although  
the conjectured rationality of the Chern-Simons invariants
would imply  $\la_{SU(3)}(X) \in \QQ$.)
Secondly,  because Casson's invariant is a finite type invariant,
one expects the same is true of the generalized Casson invariants.
This is not the case for $\la_{SU(3)};$ the computations in
\cite{bhkk} imply $\la_{SU(3)}$ is not a finite type invariant.
Nevertheless, those same computations  
support  the conjecture that $\tau$ is a finite type invariant.
Thirdly and most importantly, 
$\tau$ is   easier to compute than $\la_{SU(3)}$ 
and therefore seems more likely to satisfy a surgery formula.

To illustrate this last point, we suppose that the
moduli space $\cM$ of {\it unperturbed} flat $SU(3)$ connections on $X$ is
regular. This is equivalent to the topological assertion
$$  \text{($*$)} \quad \quad \quad\quad \quad  H^1_\al(X;su(3)) = 0
\text{ for all representations }  \al:\pi_1 X \to SU(3). \quad \quad \quad \quad \quad$$
(For example, every Brieskorn sphere of the form $\Si(2,p,q)$ satisfies
($*$).) If $X$ satisfies ($*$), then
$ \la_{SU(3)}(X) $ can be computed directly from  $\cM$
without resorting to the use of perturbations. (I.e., one can take $h=0.$)
In this case, the correction term
$\la'' _{SU(3)}$ is simply a signed sum
Atiyah-Patodi-Singer rho invariants,
but computing  $\la'' _{SU(3)}(X)$ is somewhat involved
(cf.~the computations in \cite{bhkk}).

By contrast,  whenever ($*$) holds  
the correction term $\tau''$  vanishes.
Hence,  $\tau(X)= \la'_{SU(3)}(X,0)$ 
 whenever  $\cM$ is regular.
This nice property of $\tau$ holds in the more general situation of  
homology 3-spheres satisfying
$$  \text{($**$)} \quad \quad  \quad \quad  H^1_\al(X;\CC^2) = 0
\text{ for all representations }  \al:\pi_1 X \to SU(2). \quad \quad \quad \quad \quad$$
The condition $(**)$ implies that 
$\la' _{SU(3)}(X,h)$ and $\la'' _{SU(3)}(X,h)$ 
are each independent of $h$ small.
In Lemma \ref{lemma8}, we prove that,  whenever ($**$) holds,
$\tau''(X,h)=0$ and $\tau(X) = \la'_{SU(3)}(X,h)$ 
for $h$ small.

Of course, many homology 3-spheres fail to satisfy ($**$);
examples include Brieskorn spheres of the form  $\Si(p,q,r)$ with
$p,q,r>2$. 
Thus the  integer $\la'_{SU(3)}(X,h)$ will generally depend on the choice of $h$.
Nevertheless, $\tau(X)$ agrees with $\la'_{SU(3)}(X,h)$ whenever 
the latter
is independent of $h$.  In this sense,
$\tau$ is a topological invariant of homology 3-spheres
which accounts for gauge orbits of reducible connections
only when absolutely necessary.

\bigskip \noindent
 {\bf Notation.} Throughout this paper, $X$ will be a homology 3-sphere, 
 i.e., a 
 closed,
 oriented 3-manifold satisfying $H_i(X;\ZZ) =H_i(S^3;\ZZ)$.
We denote by $\bar{X}$ the oppositely oriented homology 3-sphere.
 
Let $\la_{SU(2)}(X)$ be Casson's  original invariant, as normalized in \cite{walker}
(so $\la_{SU(2)}(X) \in  2 \ZZ$), and
$ \la_{SU(3)}(X)$  
be the  invariant defined in \cite{bh1}. 

For convenience, we use the notation for differential forms
whereby $\Om^{0+1} = \Om^0 \oplus \Om^1$,
and similarly for  
cohomology. Additionally, we use the shorthand $h^i=\dim H^i$.

Given a path $K_t,  t \in [0,1]$   of operators with discrete, real spectrum, we denote by
$SF(K_t)$ the spectral flow of $K_t$ from $t=0$ to $t=1$ using the 
$(-\ep, -\ep)$ convention. Thus, $SF(K_t)$ is the oriented intersection number
in $[0,1] \times \RR$ of the spectrum of $K_t$ with the horizontal line segment $(t,-\ep), t\in [0,1]$
for all $\ep>0$ sufficiently small.
With this convention, spectral flow is additive under composition
of paths. 

\section{Main results}
We begin with a review of 3-manifold $SU(3)$ gauge  theory.
Let $X$ be a homology sphere, $P= X \times SU(3)$, and $\th$ be the trivial
(product) connection on $P$, with covariant derivative $d$. 
We denote the space of smooth $SU(3)$ connections         
on $P$ by 
$$\cA= \{ d+ A \mid A \in \Om^1(X;su(3)) \}$$ and the gauge group of smooth bundle
automorphisms  by
$$\cG  
\cong \{g: X \to SU(3)\}.$$  This
group acts on $\cA$
by $d+A \mapsto d+gAg^{-1} + gdg^{-1}$
with quotient $\cB  = \cA / \cG$, the space of gauge orbits.  
We refer to the connection $d+A$ simply
as $A$ when no confusion can arise and
use $[A]$ to denote the $\cG$ orbit of $A\in \cA.$

Let $\cG_0$ be the identity component of $\cG$ and
set $\tB= \cA  / \cG_0.$  
Given $A \in \cA,$ we denote its $\cG_0$ orbit by
$[\![ A ] \!].$ 
Since $\cG_0$ is 
the kernel of $ \deg:\cG \to \ZZ$ (see Prop. 4.2 in \cite{bh1}),  
the natural projection $\tB \to \cB$ 
is a nontrivial connected $\ZZ$-cover (in fact the universal cover).

The Chern-Simons function $cs:\cA \to \RR$ 
is defined by the formula
$$cs(A) = \frac{1}{8 \pi^2} \int_X tr(A \wedge dA + \tfrac{2}{3} A \wedge A\wedge A).$$
Since $cs(g \cdot A)=cs (A) + \deg(g)$, the  Chern-Simons function is a well-defined
$\RR$-valued function on $\tB$, whereas on $\cB$ it takes values in
$\RR/\ZZ. $ 

The situation for the spectral flow is similar.
Choose a Riemannian metric on $X$.  Let $\cF$ denote the space of 
admissible perturbation functions (see \cite{bh1} for details).
To each pair $(A,h) \in \cA\times \cF$ we associate
 a self-adjoint, Fredholm operator
$K_{A,h}$ on  $\Om^{0+1}(X;su(3))$ defined by  the formula
\begin{equation} \label{oddsign}
K_{A,h} (\xi,a) = (d_A^* a, d_A \xi +* d_A a -4\pi \hess h(A) a).
\end{equation}
If $h=0$, then $K_{A,0} = K_A$ is simply the odd signature operator 
coupled to the connection. If $A$ is flat, then 
the Hodge and de Rham theorems identify $\ker K_A$ 
with the twisted cohomology $H^{0+1} _{\al}
(X;su(3))$, where $\al = \hol_A:\pi_1 X \to SU(3)$ is the holonomy 
representation of $A$.   Hereafter, this cohomology will be 
 denoted by $H^{0+1}_A (X; su(3))$.
More generally,  when $A$ is $h$-perturbed flat,  we can decompose $\ker K_{A,h}$
by degree and write
$$\ker K_{A,h}= H^{0} _{A}(X; su(3))\oplus
H^{1} _{A,h}(X; su(3)).$$
Note that $H^{0} _{A}(X; su(3)) = \ker (d_A:\Om^0(X;su(3)) \to \Om^1(X;su(3))$ 
is independent of $h$ and
can be identified with the Lie algebra of the stabilizer of $A$ in $\cG$.

Given  $(A_t,h_t), \ 0\leq t \leq 1,$  in $\cA \times \cF$, the spectral
flow of the path of self-adjoint operators 
$K_{A_t,h_t}, \ 0 \leq t \leq 1,$ is an integer-valued 
invariant of the homotopy class of the path 
rel endpoints.
We use the $(-\epsilon,-\epsilon)$ convention  throughout this article.
Because $\cA$ and $\cF$ are simply connected, this  spectral flow depends
only on the endpoints $(A_0,h_0)$ and $(A_1,h_1)$.

To avoid cumbersome notation, we will adopt the following
conventions.    If $A_0$  and $A_1$ are flat connections, then
$SF(A_0,A_1)$ will always mean the spectral flow of  $K_{A_t,0}$  for
$A_t, \  t\in [0,1]$ a path from $A_0$ to $A_1$.  If $A_0$ is flat (e.g.
$A_0$ is the trivial connection $\theta$) and $A_1$ is $h$-perturbed flat
for a fixed perturbation $h$, then
$SF(A_0,A_1)$ will always denote the spectral flow of $K_{A_t, h_t}$ for
$A_t,  \ t\in [0,1]$ a path of connections from $A_0$ to $A_1$ and $h_t$
a path of perturbations from $h_0=0$ to $h_1=h$, i.e. the spectral flow
of the path of self-adjoint operators from $K_{A_0,0}$ to $K_{A_1, h}$.
In all other contexts the choice of the path $h_t$ will be specified (or
obvious) and $SF(A_0,A_1)$ will denote the spectral flow of
$K_{A_t,h_t}$.

Given a connection $A$ and a gauge transformation  $g$, the index theorem
implies   that the  spectral flow of $K_{A_t,0}$  along a path $A_t$ from
$A_0=A$ to $A_1=g\cdot A$  equals $12$ deg $g$
  (for a demonstration of this, see \cite{kkr}).  Thus the function
$A\mapsto SF(\th, A)$ on connections descends to a  well-defined function
$\tB=\tB\times \{0\} \to \ZZ$ or to a function
$\cB=\cB\times\{0\}  \to \ZZ_{12}$.

Let $F_A = dA +A \wedge A$ be the curvature of the connection A and
$\cM$  the moduli space of flat $SU(3)$ connections
$$\cM = \{ A \in \cA \mid F_A  =0 \} /\cG.$$   
Note that $\cM$ is a compact subset of $\cB$
since it is homeomorphic to the space
of conjugacy classes of $SU(3)$ representations of $\pi_1 X$.
Its preimage in $\cA$ is precisely the set of critical
points of  the Chern-Simons function.
Given an admissible perturbation $h$ (see Definition 2.1 of \cite{bh1}),
a connection $A$ is called {\it $h$-perturbed flat} if it is critical point of the function
$cs+h:\cA \to \RR$. We denote
the moduli space of $h$-perturbed flat connections by $\cM_{h}$; it is   
 compact by  Lemma 8.3 of \cite{taubes}. 
 The moduli spaces of $\cG_0$ orbits of flat connections
 and $h$-perturbed flat connections are denoted $\tM$ and $\tM_{h}$.
 Since
 $\tB \to \cB$ is a $\ZZ$-cover, neither $\tM$ nor $\tM_{h}$ is compact.

Let  $\Ga_A=\{ g \in \cG \mid g \cdot A = A\}$  
be the isotropy group of  $A \in \cA$ 
and  define the subsets
$$\cA^*  = \{ A \in \cA \mid \Ga_A \cong \ZZ_3\} \quad \hbox{ and } \quad 
\cA^r= \{ A \in \cA\mid \Ga_A \cong U(1) \},$$ 
of irreducible connections and 
reducible, nonabelian connections, respectively.
Since $X$ is a homology sphere, $\cM$
decomposes as the disjoint union 
$$\cM = \{[\th]\} \cup \cM^r \cup \cM^*.$$

We denote by $\cA_{S(U(2)\times U(1))}$ the subset of $\cA$ consisting of
 connections $A$ whose holonomy reduces to the
standard $S(U(2)\times U(1))$ subgroup.
For such connections, 
$K_{A,h}$ acts diagonally with respect to the splitting
\begin{equation} \label{ohcanada}
\Om^{0+1}(X;su(3)) = \Om^{0+1} (X; \hh) \oplus \Om^{0+1} 
(X;\hh^\perp)
\end{equation}
  associated to the splitting
$su(3) = \hh \oplus \hhp$, 
where $\hh = s(u(2)\times u(1))$ and $\hhp \cong \CC^2$ is its orthogonal complement.
For any $A_0, A_1 \in \cA^r,$ choose a path in $\cA^r$ between them
($\cA^r $ is path connected) and gauge transform so the path lies in $\cA_{S(U(2)\times U(1))}.$
Then the spectral flow splits
according to the decomposition of (\ref{ohcanada})
as $$SF (A_0,A_1) =  SF_\hh(A_0,A_1) + SF_\hhp(A_0,A_1).$$
If, in addition, $A$ is $h$-perturbed flat, then  we have decompositions

\begin{equation*}
\begin{split}
H^{0}_{A}(X;su(3)) & =H^{0}_{A}(X;\hh) \oplus 
H^{0} _{A}(X; \hhp),  \\  
H^1_{A,h}(X;su(3)) &= H^1 _{A,h}(X;\hh) \oplus H^1 _{A,h}(X;\hhp).
\end{split}
\end{equation*}
(In this case, it is not hard to show that
  $H^{0} _{A}(X; \hh) \cong \RR$ and $H^{0} _{A}(X; \hhp) = 0.$)
\begin{prop} \label{prop1}
On any component $\tC \subset \tM^r,$ 
the function $\tC \to \ZZ$ defined by
$[\![ A ] \!] \mapsto SF_\hhp(\th,A)$ 
is bounded above and below. 
\end{prop}

\begin{proof} 
Since  the Chern-Simons function is constant
on components of flat connections and since 
$cs:\cB \to \RR / \ZZ$ classifies
the $\ZZ$-cover $\tB \to \cB$, 
it follows that $\tM \to \cM$ is the trivial $\ZZ$-cover.
Thus, every such $\tC$ is a homeomorphic copy  of a component of  
$\cM^r$
and is therefore compact.

Choose $[A]\in \cM^r $.  Then $A$ is nonabelian, and 
we can assume after gauge  transformation
that $A \in \cA_{S(U(2)\times U(1))}$.  
The zeroth cohomology 
$H^0 _{A}(X;su(3))$ then consists of 0-forms which are constant diagonal matrices of the form
$$\left( \begin{array}{ccc} ia  & 0&0 \\  0& ia &0  \\ 0&0 & -2ia
\end{array} \right).$$  In  particular, 
this implies that $H^0 _{A}(X;\hh^\perp)=0$ and $H^1 _A(X;\hh^\perp)$ 
is identified with the kernel of $K_{A,h}$ acting on $\hh^\perp$-valued forms.
 
 Since the dimension of the kernel of a continuous family of 
Fredholm operators is upper semicontinuous,
$h^1 _A (X; \hh^\perp )$ is a bounded function on the
compact set $\cM^r.$ 
Since we are using
the $(-\ep,-\ep)$ convention, it follows  that every
$[\![ A ] \!] \in \tC$ is contained in a neighborhood $U$  such that 
$$SF _{\hh^\perp}(\th,A) - h^1_A(X;\hh^\perp) \leq SF _{\hh^\perp}
(\th,A') \leq SF_{\hh^\perp}(\th,A)$$ for all $A' \in U$. 
Taking one such neighborhood for each $[\![ A ] \!] \in \tC$ gives an open covering of $\tC.$
Using compactness to pass to a finite subcover, we conclude that the
function $[\![ A ] \!] \mapsto SF _{\hh^\perp}(\th,A)$ is bounded above and below.
\end{proof}

For the remainder of this paper, we denote by
$C_1,  \ldots, C_n$  the connected components of $\cM^r$
and by $U_1,\ldots, U_n$ 
disjoint open sets of $\cB^r = \cA^r/\cG$ with
$C_i \subset U_i$ for $i=1,\ldots n.$
Since the Chern-Simons function $cs: \cB^r \to \RR/\ZZ$ is
 constant along the components $C_i,$  
we can choose $U_i$ small enough so that  
$ cs(U_i)$ is a proper subset of $\RR/\ZZ$.
This condition guarantees that the
restriction of  $p:\tB^r \to \cB^r$
to $p^{-1}(U_i)$ is the trivial $\ZZ$-cover.

Proposition 3.7 of \cite{bh1} insures that for $h\in \cF$
sufficiently small, we have $\cM^r _h
\subset \bigcup_{i=1}^n U_i$.
If  $[A_0], [A_1] \in U_i$, there is an
unambiguous way
to define $\hhp$ spectral flow between them by choosing gauge
representatives in
the same component of $p^{-1}(U_i)$.  (Which component of $p^{-1}(U_i)$
they lie in does not matter, since $SF_\hhp(A_0,A_1)=SF_\hhp
(g\cdot A_0, g\cdot A_1)$.)  In this situation, we will call the
representatives $A_0, A_1$ of $[A_0], [A_1]$ {\em compatible} with one
another.  Any time we refer to $\hhp$ spectral flow between
nontrivial reducible connections, we assume the connections are
compatible.

 Proposition 1 has the following consequence.

\begin{cor} \label{corollary2}
There exist $[\widehat{A}^+_i], [\widehat{A}^-_i] \in C_i, \ i=1,\dots, n,$ such
that for all $[A]\in C_i$ and any mutually compatible representatives
$A, \widehat{A}^+_i, \widehat{A}^-_i$,
\begin{equation} \label{minimouse}
\begin{split}
& SF_\hhp(\th, A)  \leq SF_\hhp (\th,\widehat{A}^+_i), \\
& SF_\hhp (\th,\widehat{A}^-_i) -h^1_{\widehat{A}^-_i}(X;\hhp) \leq SF_\hhp(\th, A) -h^1_A(X;\hhp).
\end{split}
\end{equation} \qed
\end{cor}

For generic small $h\in \cF$, $\cM_h$ is regular,
i.e., for all $[A]\in \cM_h$ the condition $H^1 _{A,h}(X;su(3))=0$ holds.
(See Section 3 of \cite{bh1} for details.)  Regularity implies that
$\cM_h$ consists of only finitely many points.

\begin{prop} \label{prop0} Suppose $h$ is a small, generic perturbation.
Define 
 \begin{equation} \nonumber
 \begin{split}
\tau''(X,h) =   \tfrac{1}{4} \sum_{i=1}^n \sum_{[A] \in \cM^r_{h} \cap U_i}  &(-1)^{SF(\th, A)}  
\left(SF_{\hhp}(\widehat{A}^+_i , A)  \right. \\
& \left. + SF_{\hhp}(\widehat{A}^-_i , A) + h^{1} _{\widehat{A}^-_i} (X;\hhp ) \right).
 \end{split}
\end{equation}  
Then $\tau''(X,h)$ is an integer and depends only on the perturbation $h$ and the
manifold $X$; in particulalr it is independent of the choice of
$ \widehat{A}^\pm_i.$
\end{prop}

\begin{proof} We first prove that $\tau''(X,h)$ depends only on the perturbation $h$.
Since we have already seen  that the spectral flow terms are gauge invariant,  
we just need to show that
$\tau''(X,h)$ is independent of the choices of $ [\widehat{A}^+_i], [\widehat{A}^-_i] \in C_i$
for $i=1, \ldots, n.$ Suppose then that
$[\widehat{B}^+_i], [\widehat{B}^-_i] \in C_i$ also 
satisfy Corollary \ref{corollary2}. Taking  lifts
$ \widehat{B}^+_i, \widehat{B}^-_i$  compatible with $ \widehat{A}^+_i, \widehat{A}^-_i,$ it follows
from additivity of spectral flow and  Corollary \ref{corollary2}
that
\begin{equation} \nonumber
\begin{split}
SF_\hhp(\widehat{A}^+_i, A) &= SF_\hhp( \widehat{B}^+_i, A), \\
SF_\hhp(\widehat{A}^-_i, A) + h^1_{\widehat{A}^-_i}(X;\hhp) 
&= SF_\hhp( \widehat{B}^-_i, A) + h^1_{\widehat{B}^+_i}(X;\hhp).
\end{split}
\end{equation} 
for all $A \in \cA^r$. This shows $\tau''(X,h)$ is independent of the choice of
$ [\widehat{A}^+_i], [\widehat{A}^-_i] \in C_i$ satisfying Corollary \ref{corollary2}.

To show  $\tau''(X,h) \in \ZZ$, we claim  that 
\begin{equation} \label{squirrel}
SF_{\hhp}(\widehat{A}^-_i , A) + SF_{\hhp}(\widehat{A}^+_i , A)   + h^{1} _{\widehat{A}^-_i} (X;\hhp ) 
\end{equation}
is divisible by 4 for all $[A] \in \cB^r.$
Additivity of the spectral flow gives that (\ref{squirrel}) equals
$$ 2SF_{\hhp}(\widehat{A}^-_i , A) - SF_{\hhp}(\widehat{A}^-_i , \widehat{A}^+_i)   + h^1 _{\widehat{A}^-_i} (X;\hhp ). $$
We claim that each of these three terms is divisible by 4.

Divisibility of the first term follows because the
$\Ga_ A \cong U(1)$ action on $\Om^{0+1}(X;\hhp)$
gives rise to a complex structure with which $K_{A_t,h_t}$ commutes
for each t, where $(A_t,h_t)$ is a path from
$(\widehat{A}^-_i,0)$ to $(A,h)$.
This implies $ SF_{\hhp}(\widehat{A}^-_i , A)$ is even.

Divisibility of the second and third terms is 
a consequence of the  following claim.
  
  \medskip
  \noindent{\it Claim.}  If $A_1$ and $A_2$ are $SU(2)\times \{ 1\}$ connections,
then $SF_{\hh^\perp} (A_1,A_2)$ and $h^1 _{A_i} (X; \hh^\perp)$ are 
divisible by four.
\medskip

To see this, identify $SU(2)$ with $Sp(1),$ the unit quaternions,
and $\hh^\perp \cong \CC^2$ with  $\HH$, the quaternions. The regular representation
of $SU(2)$ on $\CC^2$ can then be viewed    as left multiplication in $\HH$, 
and it follows that right multiplication in $\HH$
endows each eigenspace of $K_A$ with a quaternionic
 structure.  This proves the claim.
 \end{proof}

The following theorem is our main result.

\begin{thm} \label{theorem4} Suppose $h$ is a small generic perturbation.
Set
$$\tau (X)= \la'_{SU(3)}(X,h) + \tau'' (X,h),$$ where  
$\la'_{SU(3)}(X,h)$ is defined in equation (\ref{ladefn})
and  $\tau''(X,h)$ is given in Proposition \ref{prop0}.
Then $\tau(X)$ is
an integer valued  invariant
of homology 3-spheres which agrees with $\la'_{SU(3)}$ on homology 3-spheres satisfying 
{\rm ($**$)}. Furthermore
\begin{enumerate}
\item[(i)] If $\tau(X) \neq 0,$ then there exists an   
irreducible representation $\rho: \pi_1 X \to SU(3).$
\item[(ii)] $\tau(X) = \tau(\bar{X})$.
\item[(iii)]  If $X_1$ and $X_2$ are homology 3-spheres, then \\
$\tau (X_1 \# X_2) = \tau(X_1) + \tau (X_2) + 4\la_{SU(2)}(X_1) \la_{SU(2)} (X_2).$
\end{enumerate}
\end{thm}

\section{Proofs}
Although it is possible to give a direct proof of Theorem \ref{theorem4}  
based on the arguments of \cite{bh1, bh2}, 
it is in fact easier and more informative to
study the difference between $\la_{SU(3)}$ and $\tau$.
This is the approach we take in proving Theorem \ref{theorem4}.
The principal result is  Lemma \ref{lemma6}, 
where we identify  $\la_{SU(3)} - \tau$
in terms of the following more general construction.

\begin{prop} \label{prop5} 
Recall that $C_1, \ldots, C_n$ are the connected components of $\cM^r$
with disjoint neighborhoods $U_1, \ldots, U_n$ in $\cB^r$. 
Given any $\al_1, \ldots, \al_n \in \RR,$
the quantity
\begin{equation} \label{refcas}
\sum_{i=1}^n  \sum_{[A] \in \cM^r_{h} \cap U_i} 
(-1)^{SF(\th,A)} \al_i,
\end{equation}
is independent of choice of   generic small perturbation $h$.
 \end{prop}

 \begin{proof} Notice that if $\al_i=1$ for all $i$,   
the quantity (\ref{refcas}) equals  $  \la_{SU(2)}(X)$ by \cite{taubes}.
The following argument is a simple generalization of that fact.

 Suppose $\rho=h_t, 0\leq t \leq 1,$ is a generic, 1-parameter family of perturbations.
Let $$W_\rho= \bigcup_{t \in [0,1]} \cM_{h_t} \times \{t\}$$
be the parameterized moduli space. Recall that $W^r_\rho$ is a smooth
1-manifold.
If all the
perturbations $h_t$ in the path are sufficiently small,
then $W^r_\rho \subset \bigcup_{i=1}^n U_i$, hence for each $i=1, \ldots, n$,
$W^r_\rho \cap U_i$ gives a 1-dimensional cobordism from 
$\cM^r_{h_0} \cap U_i$ to
$\cM^r_{h_1} \cap U_i$ with orientations
given by the spectral flow. 
Thus each sum $  \sum_{[A] \in \cM^r_{h} \cap U_i} 
(-1)^{SF(\th,A)}$ is independent of $h$, which proves the proposition.  
 \end{proof}
 
The numbers $\al_i$ we   use to analyze the difference $\la_{SU(3)} - \tau$
are easiest to
describe in terms of the Atiyah-Patodi-Singer rho invariants.
Since $X$ is a homology 3-sphere,
every $[A] \in \cM^r$ can be represented by  
a flat $SU(2) \times \{ 1\}$ connection. 
 Given such a
connection  $A$ on $X$,
 the  rho invariant of $A$ with respect to the regular representation
 of $SU(2)$ on $\CC^2$
 can be defined by the formula (cf.~ Theorem 5.7, \cite{bhkk})
\begin{equation} \label{rhodef}
\varrho(A) = SF_{\CC^2}(\th,A) - 4 cs(A) + 2 - \tfrac{1}{2} h^1_A(X;\CC^2).
\end{equation}
Equivalently, we can replace $\CC^2$ coefficients by the subspace of $su(3)$ 
which we have been denoting by $\hh^\perp$.  
The rho invariant  $\varrho(A)$ depends only on the gauge orbit $[A]$,
not the representative.

For $i=1,\ldots, n,$ we define numbers
\begin{equation} 
\begin{split} \label{maxaldefn}
\al^+_i&= \max_{ [A]\in C_i}
\left\{\varrho(A)+ \tfrac{1}{2} h^1_{A}(X;\hh^\perp) 
\right\}, \\  
 \al^-_i&= \min _{ [A]\in C_i} 
\left\{\varrho(A)- \tfrac{1}{2} h^1_{A}(X;\hh^\perp)  \right\}.
\end{split}
\end{equation}
It is often useful to let $C_0 = \{ [\th] \}$ be the component containing 
the trivial connection and to set $\al^+_0  = 0 = \al^-_0.$

\medskip
\noindent
{\it Remark.} For each $i=1,\ldots, n,$ the connection
$[\widehat{A}^+_i]$   
 can be characterized in a gauge invariant way
as a global maximum point for the function
$C_i \to \RR$ defined by $[A] \mapsto \varrho(A) + \tfrac{1}{2}
h^1_A(X; \hhp).$
This follows by comparing equation (\ref{rhodef})
and the inequalities (\ref{minimouse}) since the 
Chern-Simons function is constant on path components 
of flat connections.
Similarly,  $[\widehat{A}^-_i]$ is a global minimum for   the function
$C_i \to \RR$ defined by
$[A] \mapsto \varrho(A) - \tfrac{1}{2} h^1_A(X; \hhp)$. 
Combining this observation with equations (\ref{rhodef}) and (\ref{maxaldefn}) shows that
\begin{eqnarray*} \label{rhomax2}
\al^+_i  
&=& SF_{\hh^\perp} (\theta, \widehat{A}^+_i) - 4 cs(\widehat{A}^+_i) +2, \\
 \al^-_i 
&=& SF_{\hh^\perp} (\theta, \widehat{A}^-_i) -4 cs(\widehat{A}^-_i) + 2 -
h^1 _{\widehat{A}^-_i} (X;\hh^\perp).
\end{eqnarray*}

 \begin{lemma} \label{lemma6}
 $\tau(X)$ is a topological invariant of homology 3-spheres.
 \end{lemma}
 \begin{proof}
Since invariance of
$\la_{SU(3)}$ is proved in \cite{bh1}, we only need to prove that
$\la_{SU(3)} - \tau$ is independent of all choices made.
Now
$$ \la'' _{SU(3)}(X,h)=\tfrac 1 2 \sum_{[A] \in \cM^r_{h}} (-1)^{SF(\th, A)} 
(SF_\hhp(\th,A)-4 cs (\hA) +2 ),$$
where $\hA$ is a reducible flat $SU(3)$ connection close to $A$.
Of course $[\hA] \in C_i$ for some $i$, and 
making compatible choices for $\widehat{A}^+_i$ and $\widehat{A}^-_i$, we see that
\begin{eqnarray*}  
\lefteqn{\la_{SU(3)} (X)     - \tau (X) =\la''_{SU(3)}(X,h) - \tau'' (X,h)} \\
&=& \tfrac 1 4 \sum_{i=1}^n \sum _{[A]\in \cM ^r _{h} \cap U_i}
(-1)^{SF(\th, A)}\left( 2SF_{\hhp}(\th, A) - 8cs (\hA)  + 4 \right. \\ 
&&\left. \hspace*{.2in} - SF_{\hhp}(\widehat{A}^+_i , A) -
SF_\hhp(\widehat{A}^-_i,A)-   h^{1}_{\widehat{A}^-_i} (X;\hhp) \right)  \\
&=& \tfrac 1 4 \sum_{i=1}^n\sum _{[A]\in \cM^r _{h} \cap U_i}
(-1)^{SF(\th, A)}\left( SF_{\hhp}(\th, \widehat{A}^+_i) -4 cs (\widehat{A}^+_i)   \right.  \\
&&\left. \hspace*{.2in}  + SF_{\hhp}(\th, \widehat{A}^-_i)
 -4cs(\widehat{A}^-_i) + 4 -   h^{1}_{\widehat{A}^-_i}(X; \hhp) \right) \\
&=&\tfrac 1 4 \sum_{i=1}^n \sum_{[A]\in \cM^r _{h} \cap U_i}
(-1)^{SF(\th, A) }\left(\al^+_i   + \al^-_i  \right).  
\end{eqnarray*}

The third step  follows by
 additivity of the spectral flow together   with the fact that 
 $cs(\hA) = cs(\widehat{A}^+_i) = cs(\widehat{A}^-_i)$, since the Chern-Simons
 function is constant along connected components of flat connections.
Now letting $\al_i=\al^+_i + \al^-_i$ and applying  Proposition \ref{prop5}
completes the proof.
\end{proof}

\begin{lemma} \label{lemma8}
 If $X$ satisfies {\rm ($**$)},
 then $\tau(X)$ equals $\la'_{SU(3)}(X,h)$ for any small generic perturbation $h$. 
 \end{lemma}
 \begin{proof}
 We show
 that $\tau''(X,h)=0$ for sufficiently small $h$ 
 whenever ($**$) holds, i.e.,
 whenever $H^1 _A (X; \hhp) = 0 $ for 
 all $[A] \in \cM^r$.
This cohomology assumption implies $H^1 _{A, h} (X;\hhp)$ also vanishes for 
every $[A]\in \cM^r _{h}$ for any  small $h.$  
(Note that the assumption of smallness of $h$ here is
stronger than the assumption needed  to define $\tau.$)
Thus $SF_\hhp(\widehat{A}^+_i, A) = 0= SF_\hhp(\widehat{A}^-_i,A)$ and $H^1_{\widehat{A}^-_i}(X;\hhp)=0$
for all $[A] \in \cM^r_{h} \cap U_i$. This shows that each summand
in the definition of $\tau''(X,h)$ vanishes for $h$ sufficiently small.
 \end{proof}

\begin{lemma}\label{lemma9}
$\tau (\bar{X}) = \tau (X)$.
\end{lemma}
 \begin{proof}
 In \cite{bh1}, it is proved that $\la_{SU(3)}(\bar{X}) = \la_{SU(3)}(X)$.
So, the lemma follows once we show that
$\la_{SU(3)}(X) - \tau(X)$ satisfies the same formula.

Reversing the orientation of $X$ changes the sign of the 
Chern-Simons function but has no effect on the perturbations.
Therefore, there is a natural correspondence between the flat moduli
spaces $\cM_h (X)$ and $\cM_{-h} (\bar{X})$.  
Obviously, if $\cM_h (X)$ is regular,  then so is $\cM _{SU(3),-h} (\bar{X})$.

The odd signature
operator $K^X_{A,h}$ acts on $\Om^{0+1}(X; su(3))$  by 
$$
K^X_{A,h}= \left[ \begin{array}{cc} 0& d_A ^*\\
d_A & *d_A - 4\pi \hess h(A) \end{array} \right].$$
Changing the orientation
of $X$ changes the sign of the Hodge star operator. 
Replacing $h$ by $-h$ as well, we see that 
 
$$ K^{\bar{X}}_{A,-h} = \left[ \begin{array}{cc} 0& d_A ^*\\
d_A & -*d _A + 4\pi \hess h(A) \end{array} \right].$$
Hence, if $ K^X_{A,h} (\xi , \eta) = \la  (\xi , \eta)$, then  
$ K^{\bar{X}}_{A,-h}(-\xi, \eta) = - \la  (-\xi , \eta).$\
Thus  switching orientations and replacing   $h$ by $-h$ 
 reflects the spectrum through zero.   
The following formula is a consequence of the $(-\ep,-\ep)$ convention:
$$SF_ {\bar{X}}(K_{A_t,-h_t}) = - SF_X(K_ {A_t,h_t}) + \dim \ker K_ {A_1,h_1} - \dim \ker K_ {A_0,h_0}.$$

Now suppose $h$ is a small perturbation and $\cM_h(X)$ is regular.
 If $[A] \in \cM^r_{h}$ then
 $H^{0}_A(X;su(3))\cong \RR$ and $H^1_{A,h}(X;su(3))  =0$ and so 
$$SF_{\bar{X}}(\th, A) = -SF_X(\th, A) + 1 -8.$$
In this formula, on   the left the spectral flow  is taken 
from $(\th,0)$  to $(A,-h)$, and on the 
right it is from $(\th,0)$  to $(A,h)$.

Further, if $\hA$ is flat and reducible, then
it is a simple exercise to prove $\varrho_{\bar{X}}(\hA) = -\varrho_X(\hA)$.
Equation (\ref{maxaldefn}) then implies that
 \begin{eqnarray*} 
 \al^+_i(\bar{X}) &=& \max_{[\hA] \in C_i} \left\{ \varrho_{\bar{X}}(\hA) +  \tfrac{1}{2} h^1_{\hA}(X;\hhp) \right\} \\
 &=& \max_{[\hA] \in C_i} \left\{ -\varrho_X(\hA) +  \tfrac{1}{2} h^1_{\hA}(X;\hhp) \right\} \\
 &=& -\min_{[\hA] \in C_i} \left\{ \varrho_X(\hA) -  \tfrac{1}{2} h^1_{\hA}(X;\hhp) \right\}  = -\al^-_i(X).
 \end{eqnarray*}
 Similarly  $\al^-_i(\bar{X}) = -\al^+_i(X)$.
 Thus
\begin{eqnarray*}  
\la_{SU(3)}(\bar{X}) - \tau(\bar{X}) &=&   \tfrac{1}{4} \sum_{i=1}^n \sum_{[A] \in \cM^r_{-h}(\bar{X}) \cap U_i} 
(-1)^{SF_{\bar{X}}(\th, A)}  
\left( \al^+_i(\bar{X})+\al^-_i(\bar{X}) \right) \\
&=&     \tfrac{1}{4} \sum_{i=1}^n \sum_{[A] \in \cM^r_{h}(X) \cap U_i} (-1)^{SF_X(\th, A)+1}  
\left( -\al^-_i(X) -\al^+_i(X)  \right) \\
&=&  \la_{SU(3)}(X) - \tau(X). 
\end{eqnarray*} 
 \end{proof}

\begin{lemma}\label{lemma10}
If $X_1$ and $X_2$ are homology 3-spheres, then\\
\hspace*{.7in}  
 $ \tau(X_1 \# X_2) =\tau(X_1)+\tau(X_2)+4\la_{SU(2)}(X_1) \la_{SU(2)}(X_2).$ 
\end{lemma}
\begin{proof}
Since $\la_{SU(3)}$ was shown to satisfy a similar 
formula in \cite{bh2}, it suffices to show additivity of
$ \la_{SU(3)} - \tau$ under connected sum. 
We first claim that the numbers $\al^+_i$ and $\al^-_i$ are additive under connected
sum. 

To make this precise, we need to set up the notation.
Set $X= X_1 \# X_2.$
Then every connection $A$ on $X$ is of the form 
 $A=A_1\#_\si A_2$, where $A_i$ is a  connection  on $X_i$ and $\si$ is
 the gluing parameter.
Furthermore, if $A$ is reducible and flat, then so are $A_1$ and $A_2$.

For $k=1,2$, let $C_0(X_k), \ldots, C_{n_k}(X_k)$
be the components of $\cM^r(X_k)$,
where $C_0(X_k) = \{ [\th_k]\}$ is the component containing the trivial connection.
The components of $\cM^r(X)$ are then given by the sets 
 $$C_{i,j}(X)= \{[A = A_1 \#_\si A_2] \mid A_1 \in C_i(X_1) \hbox{ and }
A_2 \in C_j(X_2)\}.$$ 
for $ 0 \leq i \leq n_1$ and $0 \leq j \leq n_2$.
Note that $C_{0,0}(X)$ is now the component containing the trivial connection.
For each $i,j$, we also choose an open set $U_{i,j}(X) \subset \cB^r(X)$ containing $C_{i,j}(X)$
so the collection $\{ U_{i,j}(X) \mid 0 \leq i \leq n_1, \ 0 \leq j \leq n_2\}$ is disjoint.

For $k=1,2$ and $i=1, \ldots, n_k,$ let
$\al^+_{i} (X_k)$ and $\al^-_{i} (X_k)$ be the quantities defined by equation 
(\ref{maxaldefn}).
Writing $\al^+_{i,j} (X)$ and $\al^-_{i,j} (X)$ for same numbers  
defined with respect to the components $C_{i,j}(X)$ for $X = X_1 \# X_2$, we claim 
that  
\begin{equation} \label{maxaladd}
\begin{split}
\al^+_{i,j}(X) &= \al^+_i(X_1) + \al^+_j(X_2), \\
 \al^-_{i,j}(X) &= \al^-_i(X_1) + \al^-_j(X_2).
 \end{split}
 \end{equation}
 
To see this, suppose  
$A=A_1 \#_\si A_2$ is  a reducible flat connection. Adding a 1-handle
to $\left( X_1\coprod X_2 \right) \times [0,1]$ gives a flat cobordism from 
$(X_1,A_1)\coprod (X_2,A_2)$ to $(X_1\# X_2,A)$. Since this cobordism has
no 2-handles, its signature and twisted signature vanish, and so the
Atiyah-Patodi-Singer index theorem implies that
\begin{equation} \label{rhoadd}
\varrho_{X_1\# X_2}(A)=\varrho_{X_1}(A_1)+\varrho_{X_2}(A_2).
\end{equation}
In addition, the Mayer-Vietoris principle implies
\begin{equation} \label{mv} 
H^1_A(X_1 \# X_2; \hhp) = H^1_{A_1}(X_1; \hhp) + H^1_{A_2}(X_2; \hhp). 
\end{equation}
Equations (\ref{maxaladd}) now follow by applying
 (\ref{rhoadd}) and (\ref{mv}) above to the definition (\ref{maxaldefn}).

Suppose $h_k$ is a small admissible perturbation  on $X_k$ for $k=1,2$, so that 
$\cM_{ h_k} (X_k)$ is regular.  Viewing 
$h_1$ and $h_2$ as perturbations on $X = X_1 \# X_2,$ set $h_0=h_1 + h_2$
and assume the perturbations are chosen small enough so that
$\cM^r_{h_0}(X) \subset \bigcup_{i,j} U_{i,j}(X).$

Given $[A] \in \cM^r_{h_0}(X)$, we can write $A = A_1 \#_\si A_2$
where $A_1$ is an $h_1$-perturbed flat reducible connection on $X_1$
and $A_2$ is an $h_2$-perturbed flat reducible connection on $X_2.$
Then $\cM^r _{h_0} (X)$ consists of two types of components \cite{bh2}:
 \begin{enumerate}
\item[(i)] $SO(3)$ components of the form $C=\{[A_1 \#_\si A_2]\},$
 where $[A_k] \in \cM^r _{h_k}(X_k)$ for $k=1,2$
 and $\si$ is a gluing parameter with $A_1 \#_\si A_2$ reducible. 
 \item[(ii)]  Point components of the form $C=\{[\th_1 \# A_2]\}$ or $\{[A_1 \# \th_2]$\},
 where $\th_k$ is the trivial $SU(3)$ connection over $X_k$
 and $[A_k]  \in \cM^r _{h_k}(X_k)$ for $k=1,2$.  
 \end{enumerate}
 
 Note that the intersection $\cM^r_{h_0}(X) \cap U_{i,j}(X)$ consists entirely of
 components of type (i) unless $i=0$ or $j=0$, in which case it consists of point components.
 
We first argue that  components of type (i) do not contribute 
 to $\la_{SU(3)}(X) - \tau(X)$.
 To see this, suppose $C$ is a component of type (i). Then $C \subset U_{i,j}(X)$
 for fixed $i,j>0$.
 Let $h=h_0 + tg$ be a perturbation so that    
the restriction of $g$  to $C$ is Morse.
(The existence of such functions is shown in \cite{bh2}.)  
 Then, for small $t$,  the contribution of $C$ to $\la_{SU(3)}(X) - \tau(X)$
 is given by
$$ \left(\al^+_{i,j}(X) + \al^-_{i,j}(X) \right) \sum_{p \in \crit(g|_{C})} (-1)^{\ind_p(g|_{C})},$$
which vanishes since the sum evaluates to the Euler characteristic $\chi(C)$ and
$C \cong SO(3)$. (This is similar to the proof of 
Proposition 8 in \cite{bh2}.) 
 
Thus, dropping all the $(i,j)$ terms with $i,j>0$ from the following sum
and applying equation (\ref{maxaladd})  to the remaining terms, we conclude that

\begin{equation*}
\begin{split}   \la_{su(3)}&(X) - \tau(X) \\
  & =  \tfrac{1}{4}  \sum_{i,j}
\sum_{[A]\in \cM^r _{h}(X)  \cap U_{i,j}(X)}  
(-1)^{SF_{X}(\th,A)} \left(\al^+_{i,j}(X) + \al^-_{i,j}(X)  \right)  \\
   & =  \tfrac{1}{4}  \sum_{i=1}^{n_1} \sum_{[A_1]\in \cM^r _{h_1}(X_1) \cap U_i(X_1)}  
(-1)^{SF_{X_1}(\th_1,A_1)} 
\left(\al^+_i(X_1) + \al^-_i(X_1)  \right) \\
& + \tfrac{1}{4}  \sum_{j=1}^{n_2} \sum _{[A_2]\in \cM^r _{h_2}(X_2) \cap U_j(X_2)}  
(-1)^{SF_{X_2}(\th_2,A_2)}  \left( \al^+_j(X_2)  + \al^-_j(X_2) \right) \\
& =  \la_{SU(3)}(X_1) - \tau(X_1) + \la_{SU(3)}(X_2) - \tau(X_2).
\end{split}
\end{equation*}
\end{proof}

 Since $\la_{SU(2)}(X)\in 2\ZZ$, Lemma \ref{lemma10} has the following corollary.
\begin{cor} 
The mod 16 reduction of $\tau$ is additive with 
respect to connected sum of homology 3-spheres.\qed
\end{cor}

\bigskip \noindent
{\it Proof of Theorem \ref{theorem4}.} 
Lemmas \ref{lemma6} and \ref{lemma8} show that $\tau(X)$ determines
an integer valued invariant of homology spheres $X$
which agrees with $\la'_{SU(3)}(X)$ for homology 3-spheres
satisfying ($**$). Part (i) of Theorem \ref{theorem4}
follows since $\tau(X) \neq 0$ implies
the existence of an irreducible $SU(2)$ or $SU(3)$ representation
of $\pi_1 X$, but any irreducible $SU(2)$ representation 
produces an irreducible $SU(3)$ representation via 
mapping
$SU(2) \to SO(3) \hookrightarrow SU(3).$
Part (ii) follows from
Lemma \ref{lemma9} and part (iii) from Lemma \ref{lemma10}.
\qed

Below are some computations of $\tau$ from \cite{boden2}. 
For these examples, observe that $\tau(X)$ is divisible by 2.

\begin{center}
{\sc The invariant $\tau$ for Brieskorn spheres $\Si(2,p,q)$}
\end{center}

\begin{center}
\begin{tabular}{|c||c|}  \hline
\qquad  Brieskorn sphere  \qquad & \qquad  The invariant $\tau (X)$ \qquad \\  \hline \hline
$\Si(2,3,6k\pm 1)$ & $ 3 k^2 \pm k $\\  \hline
$\Si(2,5,10k\pm 1)$ & $ 33 k^2 \pm 9k $\\ \hline
$\Si(2,5,10k\pm 3)$ & $ 33 k^2 \pm 19k+2 $\\ \hline
$\Si(2,7,14k\pm 1)$ & $ 138 k^2 \pm 26 k $\\ \hline
$\Si(2,7,14k\pm 3)$ & $ 138 k^2 \pm 62 k + 4 $\\  \hline
$\Si(2,7,14k\pm 5)$ & $ 138 k^2 \pm 102 k + 16 $\\  \hline
$\Si(2,9,18k\pm 1)$ & $ 390 k^2 \pm 58 k $\\ \hline
$\Si(2,9,18k\pm 5)$ &    $ 390 k^2 \pm 210 k + 24 $   \\  \hline
$\Si(2,9,18k\pm 7)$ &  $ 390 k^2 \pm 298 k + 52 $    \\  \hline
\end{tabular}
\end{center}

\bigskip \noindent
 {\bf Concluding remarks and open questions.}
 To better understand the relationship between $\la_{SU(3)}$ and $\tau,$
 it is helpful to compare them to
 the  $SU(2)$ invariants $\la_W$ and $\la_{BN}$
 of rational homology
 spheres defined by Kevin Walker   \cite{walker}
 and by Boyer and Nicas  \cite{bn}, respectively. 
Under suitable hypotheses, the difference $\la_W - \la_{BN}$   
can be expressed 
as a sum of  Atiyah-Patodi-Singer rho invariants 
of $U(1)$ representations \cite{CLM3}. 
A similar statement is true of $\la_{SU(3)} - \tau$, and so
it is natural to ask whether $\tau$ is the $SU(3)$ analog of the
Boyer-Nicas invariant. The answer is no, and we now explain why not.

The Boyer-Nicas approach works more generally
for   arbitrary compact Lie groups $G$.
Their idea is to define a Casson-like invariant  
by incorporating only compact components of the variety of irreducible
representations $\rho: \pi_1 X \to G$.
For that reason, no correction term is required.
In the $SU(3)$ case, their approach would yield
an integer valued invariant that presumably agrees with
$\tau$ on homology 3-spheres  satisfying ($**$).

However,  $\tau$ and the  $SU(3)$
 Boyer-Nicas invariant are not identical;
 $\tau$ involves a correction term   whereas the Boyer-Nicas invariant
does not. 
In truth, $\tau$ incorporates
 contributions from the reducible components 
 in only a very limited way. 
The proof of Lemma \ref{lemma8} shows  
that only those components of $\cM^r$ which
contain gauge orbits $[A]$ with
 $H^1_A(X;\hhp) \neq 0$    contribute to
the correction term.
These are  precisely the components  which, from
first order considerations, may contain limit points
of the irreducible stratum. The corresponding components
of the irreducible stratum would therefore  be 
excluded in the Boyer-Nicas approach
on the grounds that it is not compact.
This illustrates the fundamental difference between  
$\tau$  and the $SU(3)$  Boyer-Nicas invariant.

We conclude this paper
 with five open problems.
 \begin{enumerate}
 \item[1.] Is  $\tau(X)$ divisible by 2 for all homology 3-spheres?  
If not, is $\tau$ mod $2$ a homology cobordism invariant? 
 \item[2.]  Is $\tau$  a finite-type invariant?
 \item[3.]  Find a Dehn surgery formula for  $\tau$. 
 \item[4.] Compute $\tau$ for homology spheres that
 do not satisfy ($**$), e.g., Brieskorn spheres $\Si(p,q,r)$ with $p,q,r>2.$  
\item[5.] Develop an $SU(3)$ Floer theory and
relate its  Euler characteristic to $\tau.$
 \end{enumerate}




\end{document}